# Algorithmic Monotone Multiscale Finite Volume Methods for Porous Media Flow

Omar Chaabi and Mohammed Al Kobaisi *

*Petroleum Engineering Department, Khalifa University of Science and Technology, Abu Dhabi 127788, UAE*

## Abstract

Multiscale finite volume methods are known to produce reduced systems with multipoint stencils which, in turn, could give non-monotone and out-of-bound solutions. We propose a novel solution to the monotonicity issue of multiscale methods. The proposed algorithmic monotone (AM-MsFV/MsRSB) framework is based on an algebraic modification to the original MsFV/MsRSB coarse-scale stencil. The AM-MsFV/MsRSB guarantees monotonic and within bound solutions without compromising accuracy for various coarsening ratios; hence, it effectively addresses the challenge of multiscale methods' sensitivity to coarse grid partitioning choices. Moreover, by preserving the near null space of the original operator, the AM-MsRSB showed promising performance when integrated in iterative formulations using both the control volume and the Galerkin-type restriction operators. We also propose a new approach to enhance the performance of MsRSB for MPFA discretized systems, particularly targeting the construction of the prolongation operator. Results show the potential of our approach in terms of accuracy of the computed basis functions and the overall convergence behavior of the multiscale solver while ensuring a monotone solution at all times.

## 1. Introduction

The problem of multiphase flow in porous media is considered multiscale in nature due to the physical processes occurring over a wide range of spatio-temporal scales. Using conventional simulation approaches, it might be impractical to resolve all these scales with a single high-resolution grid. This prompted the development of numerous multiscale methods over the past 50 years. The idea pertains to using multiple computing grids at different scale levels to arrive at the solution upon the desired or finest scale. The efficiency of employing such approaches was examined through numerous numerical experiments early on [1]. An operator splitting on multiple grids (OSMG) technique [2] was developed for simulations in heterogeneous reservoirs. A finite element, with Hermite cubic basis functions, was used as the coarse grid pressure solver. The solution is then interpolated using splines to the fine grid where the transport problem was solved. Two time scales were employed, a large pressure solve time step and a small time step for carrying out the saturation solution. The procedure did not guarantee local mass conservation. A dual-mesh method for a 2-D single phase model was presented in [3]. The pressure field was first computed on the coarse grid, fine grid velocities were then obtained by solving local flow problems with approximate pressure boundary conditions obtained from the coarse solution. Later on, the method was extended to water-oil flow [4]. A similar approach that takes into account the effects of gravity was presented in [5]. The pressures on the coarse grid along with the saturation solution on the fine grid were obtained concurrently with a single time step. A nested gridding approach that computes approximate pressures and uses streamline based simulation to advance the saturations was put forward in [6]. In another study [7], a subgrid approach employing numerical Green's functions and mixed finite element methods, was used to upscale transmissivities and achieve local mass conservation. An extension of the dual-mesh technique to three dimensions incorporating gravity effects and wells was introduced in [8]. Improvements were noted based on the upscaling technique and the flux redistribution using transmissibility weighting schemes.

* Corresponding author at: PE Department, Khalifa University of Science and Technology, Abu Dhabi 127788, UAE
  *E-mail address*: mohammed.alkobaisi@ku.ac.ae (M. Al Kobaisi)

In 1997, Hou and Wu [9] presented a new paradigm to multiscale methods. The authors developed a finite element method (MsFE) whereby numerical bases functions are constructed which can capture the small-scale details effectively. However, the method did not ensure locally conservative reconstruction of the fine-grid velocity field. A similar approach was later devised using a multiscale mixed finite element method (MsMFE) in which locally conservative reconstruction of the fine-grid velocity field was demonstrated [10]. A variant of the MsMFE method, specifically focusing on blocks containing sources, was provided in [11]. Extension of the MsMFE method to non-uniform coarse grids and to corner-point and other unstructured grids can be found in references [12] and [13].

Jenny and Tchelepi [14], introduced a multiscale finite-volume method (MsFV) which proved to be seminal for modern day multiscale reservoir simulation. The MsFV method is based on a control volume finite difference approach in which dual grids are used to construct the multiscale basis functions. It is conservative by construction and closer to reservoir simulation practitioners since it is finite-volume based. Herein, we mainly focus on the finite-volume type methods; specifically, we concentrate on the original Multiscale Finite Volume (MsFV) method [14]–[19], and the state-of-the-art Multiscale Restricted Smoothed Basis (MsRSB) method [20]. Generally, multiscale methods use numerically computed basis functions, obtained by solving localized flow problems, to effectively integrate unresolved subscale effects into macroscale models while maintaining consistency with the underlying flow equations. Initially, multiscale methods were introduced as a single-pass alternative to upscaling [14]. Later on, they were incorporated into iterative schemes aimed at eliminating non-physical oscillations [21], [22]. Moreover, multiscale methods were combined with fine-scale smoothers in iterative procedures, resembling a two-level algebraic multigrid method [18], [23], [24]. Multiscale methods have now advanced to a point where they are used commercially to speed up the simulation of large and complex geomodels [25]–[27].

The original MsFV and the state-of-the-art MsRSB provide approximate pressure solutions that are generally in good agreement with reference fine-scale solutions across various heterogeneous test cases. However, multiscale methods are inherently multipoint due to the expansion of support regions in MsRSB (or dual grids in MsFV) beyond the boundaries of the coarse control volume [28]. As a result, multiscale methods can, and often will, produce pressure solutions with non-physical peaks which, in turn, affect the quality of the reconstructed velocity fields [29], [30]. This issue is further exacerbated when the underlying fine-scale contains regions with high permeability contrasts, the fine scale meshes are with high aspect ratios, or a combination thereof. To eliminate the non-physical oscillations, Hajibeyji et al., [31], [32] proposed using relaxation to smooth out the high frequency errors, thus, enhancing the quality of the reconstructed fine-scale velocity fields. Other approaches focused on localized modifications to the coarse-scale stencil to overcome the lack of monotonicity [33]–[35].

In this work, we tackle the lack of monotonicity of multiscale methods in a generic and systematic way. We propose an algorithmic monotone MsFV/MsRSB operator that guarantees strictly monotone and within bound pressure solutions without compromising the accuracy of the solution. The performance of the algorithmic monotone operator in iterative setups is also investigated and compared to MsRSB. Furthermore, we extend the algorithm to include MPFA formulations and show how the algorithm can enhance the performance of MsRSB for such consistent discretized systems.

The remainder of the paper is structured as follows. **Section 2** introduces the mathematical models and discusses the monotonicity issue of multiscale methods. In **Section 3**, we introduce the AM-MsFV/MsRSB framework and the underlying algorithmic choices. Then, in **Section 4**, we test the performance of the AM-MsFV/MsRSB and report the numerical results including comparisons to original methods and published results. In **Section 5**, we extend the AM-MsRSB to MPFA discretized systems and we introduce a new modification to the prolongation construction step. Finally, **Section 6** presents the conclusions of this work.

## 2. Model Problems

Multiscale methods are primarily designed to compute approximate solutions of second-order elliptic equations with heterogeneous coefficients. To clearly show the ideas of multiscale methods, it is sufficient to work with the single phase pressure equation under simplified physics. Assuming incompressible flow and ignoring the gravitational effects, the pressure equation is equivalent to the variable coefficient Poisson equation and can be written as,

$$\nabla \cdot \big( \mathbf{K}(x)\nabla p(x) \big) = q(x), \quad x \in \mathbb{R}^d, \quad \mathbf{K}(x) \in \mathbb{R}^d \times \mathbb{R}^d. \tag{1}$$

Here, $p$ is the pressure, and $q$ denotes the source term. The permeability tensor **K**, often exhibiting a complex spatial structure, can contain discontinuities of varying orders of magnitude and patterns spanning numerous length scales. Thus, introducing the multiscale nature into the problem. For a computational domain $\Omega \in \mathbb{R}^d$ that is d-dimensional, problem **(1)** is considered well-posed, provided that appropriate boundary conditions are imposed at $\partial\Omega \in \mathbb{R}^{d-1}$. The pressure equation is discretized using the de facto standard two-point flux approximation (TPFA), which results in a discrete equation for each cell

$$\sum_j T_{ij}(p_i - p_j) = q_i, \tag{2}$$

where $T_{ij}$ is the transmissibility, which approximates the flux across the shared interface between two neighboring computational nodes $i$ and $j$. In a two-dimensional setup of Cartesian grids, TPFA gives rise to the standard five-point finite-difference stencil. Assembling the discrete equations **(2)** for each cell results in the following linear system

$$\boldsymbol{A}\,\boldsymbol{p} = \boldsymbol{q}, \tag{3}$$

where $\boldsymbol{A} \in \mathbb{R}^{d \times d}$ is the coefficient matrix, $\boldsymbol{q} \in \mathbb{R}^d$ is a vector of source terms, and $\boldsymbol{p} \in \mathbb{R}^d$ is the vector of cell-averaged pressures that we wish to solve for.

TPFA is only consistent for **K**-orthogonal grids, i.e., grids which are orthogonal with respect to the permeability tensor **K**. Thus, other schemes, such as the multipoint flux approximation (MPFA) and MPFA-like schemes [36]–[40], nonlinear TPFA and nonlinear MPFA schemes [41]–[44], are required to capture tensorial permeability flow and handle non-**K**-orthogonal grids. In MPFA schemes, as the name implies, fluxes are approximated with multiple points compared to the two-point TPFA. This would ultimately lead to a denser coefficient matrix **A**, but would overcome the issues of **K**-orthogonality to a certain degree.

### 2.1. Multiscale formulation

In terms of generating coarse grids, both the original MsFV and the state-of-the-art MsRSB inherently adhere to similar principles, although they might employ different preprocessing algorithms. The primary distinction between the two schemes lies in the computation of the basis functions, i.e., the prolongation operator which is discussed below. If we define the discrete fine-scale domain as the collection of cells $\{\Omega_i^f\}_{i=1}^n$, and define a coarser grid based on a defined partition $\{\Omega_j^c\}_{j=1}^m$, such that each fine cell in $\Omega^f$ belongs to one cell in $\Omega^c$. We can solve **(3)** with reduced degrees of freedom on the coarse grid by defining $\boldsymbol{p}_c$ as coarse pressure values in $\Omega^c$. Moreover, to allow mapping between $\Omega^f$ and $\Omega^c$, we introduce two transfer operators, namely the prolongation ($\boldsymbol{P}_f^c$) and the restriction ($\boldsymbol{R}_c^f$). $\boldsymbol{P}_f^c$ is a $n \times m$ matrix that maps from coarse-scale to fine-scale $\{\Omega^c\} \to \{\Omega^f\}$, whereas $\boldsymbol{R}_c^f$ is a $m \times n$ matrix that maps in the opposite direction $\{\Omega^f\} \to \{\Omega^c\}$. If we assume for now that $\boldsymbol{P}_f^c$ and $\boldsymbol{R}_c^f$ are readily available, we can find an approximate fine-scale pressure $\boldsymbol{p}_f$ on $\Omega^f$ from any pressure distribution on $\Omega^c$,

$$p_f = P_f^c p_c. \tag{4}$$

Now if we use the approximate $p_f$ instead of $p$ in (3) and use $R_c^f$ to reduce the number of equations so that it matches the degrees of freedom we obtain

$$R_c^f \left( A_f \left( P_f^c p_c \right) \right) = (R_c^f A_f P_f^c) p_c = A_c p_c = R_c^f q_f = q_c. \tag{5}$$

The quality and physical interpretation of this system is influenced by the choice of the prolongation and restriction operators.

### 2.1.1. Prolongation and restriction operators

The prolongation operator can be looked at as a concatenation of basis functions, in which a column $j$ in $P_f^c$ entails the basis functions corresponding to $\Omega_j^c$. The set of basis functions for a coarse block is simply the local pressure response in the block when a unit pressure is imposed in the block itself. The original MsFV uses a dual coarse grid and solves reduced local problems to compute basis functions. On the other hand, MsRSB initiates $P_f^c$ as the characteristic function for each coarse block, then uses an algebraic smoother to transform $P_f^c$ so that its columns adjust to the discretized flow equations. Furthermore, a special localization strategy is applied to avoid the stencil growth outside of the predefined support region, and this adjustment acts as a sort of reduced boundary condition [45]. We refer the interested reader to the following references for more details on how to obtain the prolongation operator for each scheme [14], [20], [28].

On the other hand, there are two choices for the restriction operator that are frequently used in the literature;

$$(R_{CV})_{ij} = \begin{cases} 1, & \text{if } \Omega_i^f \in \Omega_j^c, \\ 0, & \text{otherwise,} \end{cases} \quad \text{or} \quad R_G = P^T.$$

The first choice, which is the characteristic function of each coarse control volume, is known as the control-volume restriction operator and it was introduced in the original MsFV formulation [14]. Using the control-volume restriction operator gives a coarse system matrix where the transmissibility term between two neighboring coarse blocks, say $i$ and $j$, is the sum of the estimated fluxes from the basis functions over the coarse interface $\Gamma_{ij}^c$. The use of $R_{CV}$ is a requirement for the flux reconstruction step, which is needed to find a mass conservative velocity field at the fine-scale [46]. However, using $R_{CV}$ can lead to non-physical pressure solutions at the coarse-scale leading to severe errors. For highly heterogeneous porous media, using $R_{CV}$ give pressure solutions that are non-monotone and violate the maximum principle [19]. Therefore, when pressure stability is a concern, $R_G$ is usually favored [18]. The second choice, $R_G$, constructs the restriction operator using the Galerkin principle, which is the typical choice for multigrid methods.

**2.2. Iterative multiscale**

Iterative multiscale solvers have been devised and applied to either systematically lower the residual towards zero or to eliminate undesired values and non-monotonicity from the solution [18], [21], [23]. Setting up an iterative multiscale solver is similar to deriving a two-level algebraic multigrid method. Let the solution at step $k$ of the iterative scheme be $\boldsymbol{p}^k$ and define the residual $\boldsymbol{r}$ as:

$$\boldsymbol{r}^k = \boldsymbol{q} - \boldsymbol{A}\boldsymbol{p}^k. \tag{6}$$

Additionally, assume that we have a function $S(\boldsymbol{A}, \boldsymbol{b})$ that executes one or more smoothing iterations, at the fine-scale, for a given matrix A and right-hand-side b. Note that any affordable iterative solver that effectively eliminates high frequency errors from the solution is referred to as a smoother. Considering this, we may define a two-step preconditioner that first removes high frequency errors using the smoother and then addresses low frequency errors by a multiscale coarse grid correction,

$$\boldsymbol{p}^{k+1/2} = \boldsymbol{p}^k + S(\boldsymbol{A}, \boldsymbol{q} - \boldsymbol{A}\boldsymbol{p}^k) = \boldsymbol{p}^k + S(\boldsymbol{A}, \boldsymbol{r}^k), \tag{7}$$

$$\boldsymbol{p}^{k+1} = \boldsymbol{p}^{k+1/2} + \boldsymbol{P}_f^c \boldsymbol{A}_c^{-1} \boldsymbol{R}_c^f (\boldsymbol{q} - \boldsymbol{A}\boldsymbol{p}^{k+1/2}) = \boldsymbol{p}^{k+1/2} + \boldsymbol{P}_f^c \boldsymbol{A}_c^{-1} \boldsymbol{R}_c^f \boldsymbol{r}^{k+1/2}. \tag{8}$$

Once the coarse grid correction has been computed, equation **(8)**, the solution is mass conservative over the coarse domain, given that $\boldsymbol{R}_{CV}$ is used. However, and as proven in [23], it is sufficient to only conclude the iterations by $\boldsymbol{R}_{CV}$ to guarantee mass conservation at the coarse-scale. The iterative multiscale solver can be used to obtain solutions with any accuracy. If the linear tolerance is very tight, it is also common to employ Krylov-based algorithms to speed up the solution process [23].

**2.3. Multiscale lack of monotonicity**

To better understand the source of non-monotonicity, we could envisage the multiscale method's reduced system as a single-scale multipoint control-volume system. Multiscale methods are inherently multipoint due to the expansion of support regions in MsRSB (or dual grids in MsFV) beyond the confines of the coarse control volume [28]. This means that the computed basis functions overlap with multiple neighboring coarse blocks, which also introduces effective permeabilities at the coarse grid that are not aligned with the grid [34]. Furthermore, multipoint schemes are generally not monotone [47]. Accordingly, the reasoning employed in single-scale methods can be extended to multiscale methods and informally state that: there is no guarantee that multiscale methods provide solutions that are monotone and obey the maximum principle.

Attempts to address the lack of monotonicity of the MsFV operator resulted in the development of iterative formulations. While other studies focused on locally modifying the coarse-scale stencil to get monotone multiscale solutions [33], [35]. Wang et al., [35] proposed a monotone MsFV (m-MsFV) based on a local stencil-fix approach. In their approach, critical coarse-scale connections are detected and replaced by TPFA-based coarse-scale transmissibilities. Furthermore, the authors nicely highlight and identify the causes of the non-physical oscillatory behavior of the MsFV operator. The large unphysical pressure peaks in the original MsFV method can be attributed to the localization assumptions. The oscillatory behavior becomes more evident in regions with high contrasts in the permeability between adjacent cells [28]. While the oscillatory behavior is more pronounced in the MsFV operator, the state-of-the-art MsRSB operator also does not guarantee monotone and within bounds solutions. This limitation arises from the same reasons mentioned earlier, which apply to multiscale methods in general.

In the coarse-scale system $\boldsymbol{A}_c$, an entry $a_{ij}^c$ represents a sum of the estimated fluxes from the basis functions over the coarse interface $\Gamma_{ij}^c$. For cases when the coarse node lies in a low permeability region, the net incoming flux at $\Gamma_{ij}^c$ from cell $j$ could be negative, which results in a positive off-diagonal value for $a_{ij}^c$. Moreover, the net outgoing fluxes of cell $i$ within its own control volume could be relatively small, leading to a lack of diagonal dominance in the corresponding row of the coarse-scale system. Therefore, it is possible to algebraically detect the specific rows and entries responsible for the non-physical oscillatory behavior of the multiscale coarse solution. Given that the basis functions are typically monotone with values between 0 and 1, monotonicity violations at the fine-scale can only be caused by violations at the coarse-scale [35]. In other words, fixing the monotonicity issue at the coarse-scale will guarantee monotone approximate solutions at the fine-scale. In the following section, we propose a remedy to the lack of monotonicity of multiscale methods.

## 3. Algorithmic monotone multiscale methods

Building on the findings from the previous section, it is evident that the non-monotone behavior of multiscale methods is attributed to the existence of positive off-diagonal entries and the lack of diagonal dominance. Both of which are deviations from M-matrix properties, and they collectively contribute to the non-monotone behavior. Accordingly, we aim to enforce M-matrix properties at the coarse-scale to obtain monotone coarse-scale solutions. This manipulation is generally acceptable in a multiscale setup because the primary goal is to obtain excellent approximate solutions, rather than exact solutions. However, improvident manipulation might impact the ability of multiscale methods to generate satisfactory approximate solutions.

Next, we explain the mathematical motivation behind our algorithmic choices for carefully enforcing M-matrix properties at the coarse-scale. Although we show the applicability of our algorithm for both MsFV and MsRSB, our algorithmic choices primarily target the MsRSB formulation which obtains the prolongation operator in a similar fashion to smoothed-aggregation algebraic multigrid (SA-AMG) methods [48]. Moreover, our work was inspired by parallel multigrid applications where algebraic stencil collapsing is performed to obtain non-Galerkin coarse operators [49], [50].

In algebraic multigrid (AMG), which shares some similarities with multiscale methods, a primary objective is to ensure that the prolongation (i.e., interpolator) operator contains the near null space components within its range space [51]. This, in turn, provides the expectation that the coarse-scale system produced by AMG can accurately capture the near null space components of the fine-scale system, i.e., vectors $x$ such that $Ax \approx 0$ [52]. Such vectors coincide with the error components that are slow to resolve using a relaxation, and are known to be closely related to the algebraically smooth modes [53], [54]. Informally, multiscale coarse approximations can be thought of as the equivalent of the coarse grid correction in multigrid.

Let $A_c$ be our coarse-scale system

$$A_c = R_c^f A_f P_f^c, \tag{9}$$

where $P_f^c$ is the restricted smoothed basis of MsRSB and the restriction is either $R_{CV}$ or $R_G$. Moreover, assume that the near null space is accurately represented in the range space of $P_f^c$. Let us also assign $A_m$ as a near null space preserving and monotone coarse-scale operator, i.e., with enforced M-matrix properties where the positive off-diagonal entries are eliminated (or weakened). Let $v_c \in \mathbb{R}^m$ be a coarse-scale vector, then the goal is to obtain $A_m$ such that $A_m v_c \approx A_c v_c$. Generally, the equality will not hold for all $v_c$ because $A_m$ has fewer off-diagonal entries due to the elimination of positive entries. However, it will be enough for the equality to hold when $v_c$ corresponds to a smooth mode.

Therefore, the aim is to obtain $A_m$, while enforcing M-matrix properties and preserving the near null space, such that it scales a vector of smooth modes just like $A_c$ would. In the context of scalar PDEs, the near null space or the low energy eigenmodes are often represented by a constant vector [55]. In this study, the near null space is simply a vector of ones, $\mathbf{1} \in \mathbb{R}^m$. Subsequently, we aim to enforce the following equalities:

$$A_m \mathbf{1} = A_c \mathbf{1} \text{ and } (A_m)^T \mathbf{1} = (A_c)^T \mathbf{1}. \tag{10}$$

This can be thought of as a way of forcing the null spaces of both operators to coincide, and this is done by modifying the diagonal entry to account for the dropped positive off-diagonal entries. To achieve that, we start by letting $A_m$ be a perturbation of $A_c$,

$$A_m = A_c + B. \tag{11}$$

Then, we construct $B$ such that it satisfies $B1 = 0$ and $(B)^T 1 = 0$. For instance, say that we want to eliminate a positive off-diagonal entry $a_{ij}^c$, then we add to $A_c$ a matrix $B$ of the following form,

$$B = \begin{pmatrix} wa_{i,j}^c & -wa_{i,j}^c \\ -wa_{i,j}^c & wa_{i,j}^c \end{pmatrix} \begin{matrix} i \\ j \end{matrix} \tag{12}$$

where $w$ is a user defined weight. When $R_G$ is used, we end up with a symmetric $A_c$ and the choice of $w = 1$ will result in total elimination of the positive off-diagonal entries. While complete elimination of positive off-diagonal entries ensures monotone solutions, it is sometimes sufficient to only weaken these entries by choosing $w < 1$ or even disregarding small positive off-diagonals [35]. To address this, we introduce an indicator to identify problematic positive off-diagonal entries. This indicator, $\zeta_{ij}$, is defined as the ratio between the entry's value and the corresponding diagonal term ($\zeta_{ij} = a_{ij}^c / a_{ii}^c$). Additionally, a user-specified threshold, $\varepsilon$, is defined, and all entries satisfying $\zeta_{ij} > \varepsilon$ are considered problematic and must be addressed by the algorithm. The steps for constructing the proposed algorithmic monotone multiscale operator are summarized in **Algorithm 1** below.

---

**Algorithm 1** Algorithmic monotone MsFV/MsRSB (AM-MsFV/MsRSB)

1: Set up a course grid and support regions
2: Compute basis functions
3: Choose Restriction operator  ⇒  $R_G$ for iterative applications and $R_{CV}$ for conservative approximations
4: Construct the coarse-scale system  ⇒  Eq. (5)
5: Locate all positive off-diagonal entries, $Pos\_entry$
6: **for** $i = 1: numel(Pos\_entry)$ **do**
7:    **if** $\zeta_{ij} > \varepsilon$ **then**
8:       Specify a weight value $w$
9:       Construct the perturbation matrix $B$ as follows:  ⇒  Eq. (12)
10:       $b_{ij} \leftarrow b_{ij} - w * a_{ij}^c$
11:       $b_{ii} \leftarrow b_{ii} + w * a_{ij}^c$
12:       $b_{ji} \leftarrow b_{ji} - w * a_{ij}^c$
13:       $b_{jj} \leftarrow b_{jj} + w * a_{ij}^c$
14:    **end if**
15: **end for**
16: $A_m = A_c + B$  ⇒  Eq. (11)

Multiscale methods can be effectively used in three different modes [28]: (a) as a fine-scale linear solver, where the residual is reduced to very tight tolerance; (b) as an approximate solver that only reduces the residual to a relaxed prescribed tolerance but still guarantees a mass-conservative approximation; or (c) as a one-step approximate solver that gives conservative fluxes on the fine-scale, with no control on fine-scale residual. Note that applications (a) and (b) utilize multiscale methods in an iterative formulation. Towards that end, our algorithm ensures that the monotone $A_m$ yields a similar effect to $A_c$ with respect to algebraically smooth errors not reduced by relaxation, and in turn, ensures that it does not degrade the efficiency of the solver when using $A_m$ instead of $A_c$. Regarding application (c), our algorithm enhances the accuracy of pressure approximations as it guarantees monotone pressure solutions that abide by the maximum principle. Moreover, when employed with $R_{CV}$, the addition in equation **(11)** maintains mass balance on the coarse grid, allowing for conservative flux computation at the coarse-scale. Furthermore, and as in the original formulation, we can also construct conservative fluxes on the fine-scale by solving additional local problems where the conservative coarse-scale fluxes are imposed as Neumann boundary conditions along the coarse boundaries.

## 4. Numerical Experiments

We implemented our algorithm using Matlab Reservoir Simulation Toolbox [28], [47], [56], and leveraged its MsRSB implementation, plotting functions, and various other helpful utilities. In this section, we present the results of several test cases to validate our algorithm and demonstrate the potential of the proposed AM-MsFV/MsRSB. To evaluate the spatial accuracy of AM-MsFV/MsRSB, we examine the difference between the multiscale approximations and the fine-scale reference solution. This is done by calculating the scaled $L^2$ and $L^\infty$ norms,

$$\|p^f - p^{ms}\|_2 = \sqrt{\frac{\sum_{i=1}^n |p_i^f - p_i^{ms}|^2}{\sum_{i=1}^n |p_i^f|^2}}, \qquad \|p^f - p^{ms}\|_\infty = \frac{\max |p_i^f - p_i^{ms}|}{\max |p_i^f|}, \qquad (11)$$

where $p_i^f$ is the fine-scale reference pressure and $p_i^{ms}$ is the multiscale approximate pressure prolonged to $\Omega^f$.

### 4.1. Multiscale methods as one-step approximate solvers

In this section, we test the performance of our algorithm on multiscale methods when used as a one-step alternative to upscaling, i.e., without employing subsequent iteration cycles. Since our interest is not on investigating the iterative performance, it would make sense to use the control volume restriction operator, $R_{CV}$. Our algorithm was designed to target the lack of monotonicity in the state-of-the-art MsRSB method. However, for completeness, we first apply the algorithm to the original MsFV operator.

### 4.1.1. Test case 1: SPE 10 bottom layer – MsFV

In this example, we consider flow in the bottom layer of the 10th SPE comparative solution project [57]. The horizontal layer has (60 × 220) fine cells. We use a coarsening ratio of (5 × 10) to partition the fine-scale grid, giving rise to a (12 × 22) coarse-scale system. We impose unit Dirichlet boundary conditions to drive the flow from left to right. The permeability distribution and the reference fine-scale pressure solution are depicted in **Figure 1**.

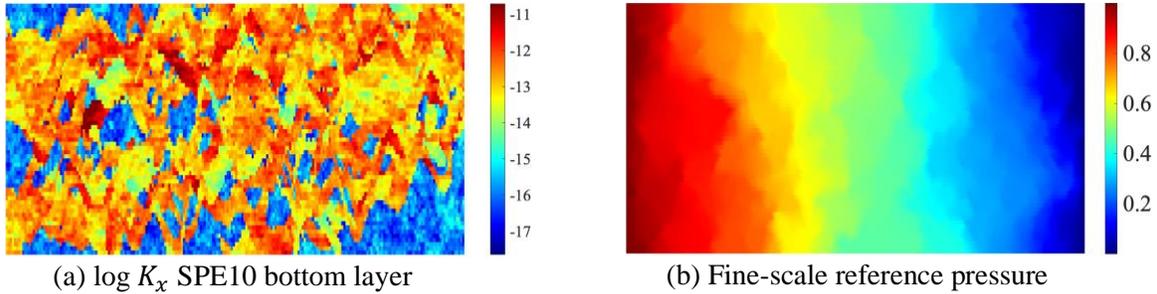

(a) log $K_x$ SPE10 bottom layer          (b) Fine-scale reference pressure

**Figure 1.** Permeability and reference pressure solution for the SPE 10 bottom layer (layer 85). A unit pressure is imposed from the left boundary to the right boundary.

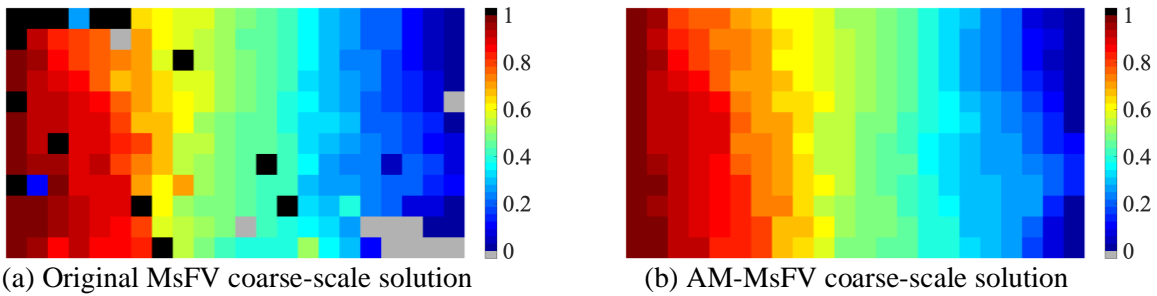

(a) Original MsFV coarse-scale solution          (b) AM-MsFV coarse-scale solution

**Figure 2.** Comparison of the coarse-scale pressure solutions computed by the original MsFV and the AM-MsFV.

Next, we use MsFV and algorithmic monotone MsFV (AM-MsFV) to compute approximate pressure solutions. **Figure 2** shows the coarse-scale solutions of both methods, with the colorbar adjusted to display values out-of-bound with distinct colors. Specifically, cells with pressures above one are shown in black and the ones with values below zero are shown in silver. Since we are solving an elliptic problem, solutions should obey the maximum principle i.e., they should not go out-of-bound [0,1]. We can see that the original MsFV exceeds those bounds at several locations which also reflects on the quality of the prolongated solutions, see **Figure 3**. Meanwhile, a monotone pressure solution is obtained by AM-MsFV at the coarse-scale with $\varepsilon = 0.001$ and $w = 1.5$. Furthermore, **Figure 4** shows the pressure surface plots for the reference pressure solution, the original MsFV, and the AM-MsFV. The surface plots highlight how our algorithm can eliminate the non-physical oscillations resulting from the localization assumptions of the MsFV operator. Quantitatively, the AM-MsFV shows $L^2$ and $L^\infty$ norms of 0.033 and 0.13, respectively and thus, emphasizing the operator's ability to provide monotone, within bounds, and accurate pressure approximations.

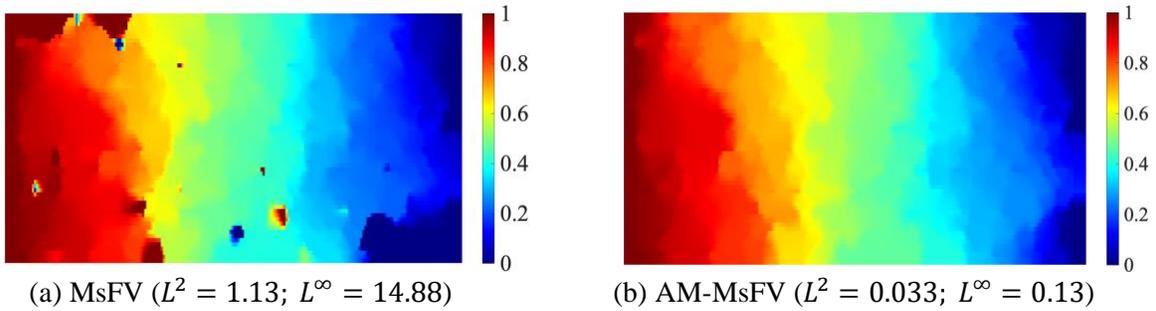

(a) MsFV ($L^2 = 1.13$; $L^\infty = 14.88$)　　　　(b) AM-MsFV ($L^2 = 0.033$; $L^\infty = 0.13$)

**Figure 3.** Comparison of the prolongated pressure solutions computed by the original MsFV and the AM-MsFV.

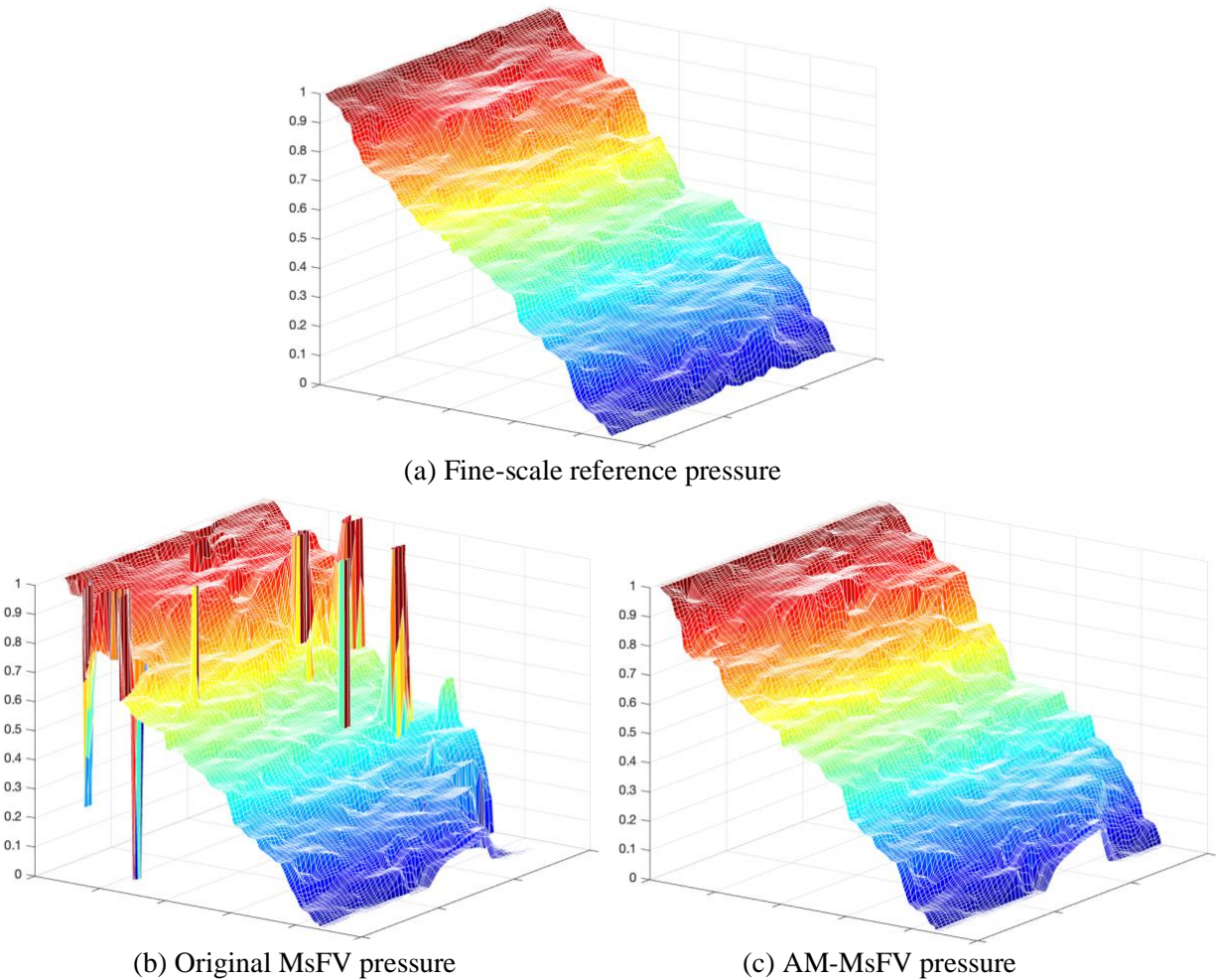

(a) Fine-scale reference pressure

(b) Original MsFV pressure　　　　(c) AM-MsFV pressure

**Figure 4.** Pressure surface plots for (a) fine-scale reference; (b) original MsFV; (c) AM-MsFV.

### 4.1.2. Test case 2: SPE 10 bottom layer – MsRSB

In this example, we once again investigate the bottom layer of SPE 10, but this time we solve using both the original MsRSB and the AM-MsRSB. With the coarsening ratio utilized in case 1 (5 × 10), the original MsRSB only exhibits slight off-bounds behavior near the inlet and outlet boundaries. Therefore, for this test case we opt for a coarsening ratio of (3 × 5) which leads the MsRSB operator to produce solutions with non-physical oscillations in the middle of the domain. The coarse-scale and prolongated pressure solutions are shown in **Figures 5** and **6**.

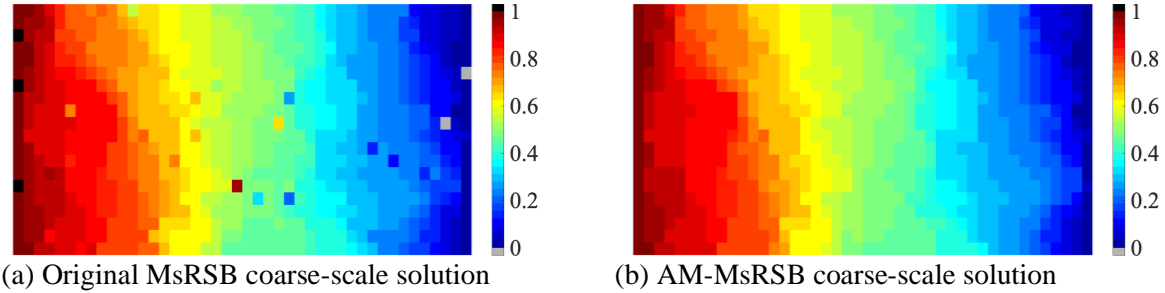

(a) Original MsRSB coarse-scale solution        (b) AM-MsRSB coarse-scale solution

**Figure 5.** Coarse-scale pressure solutions computed by the original MsRSB and the AM-MsRSB.

Pressure solutions computed by the original MsRSB method contains local oscillations and few out-of-bound values. In contrast, the AM-MsRSB operator produces a strictly monotone solution without violating the maximum principle when the threshold and weight are set to $\varepsilon = 0.1$ and $w = 1$, respectively. As evident from the pressure errors ($L^2 = 0.0093$; $L^\infty = 0.045$), the AM-MsRSB method can provide a monotone pressure solution without compromising on accuracy. **Figure 7** shows the surface plots and we can see that the AM-MsRSB provides within bounds solutions with no non-physical local peaks.

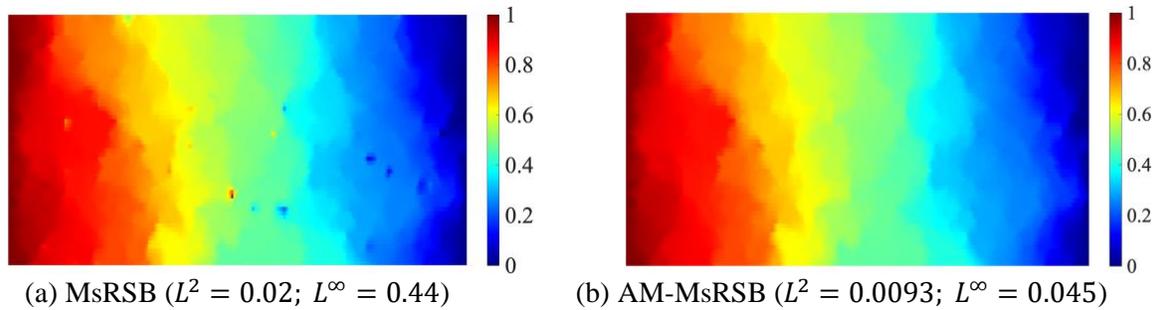

(a) MsRSB ($L^2 = 0.02$; $L^\infty = 0.44$)        (b) AM-MsRSB ($L^2 = 0.0093$; $L^\infty = 0.045$)

**Figure 6.** Comparison of the prolongated pressure solutions computed by the original MsRSB and the AM-MsRSB, refer to Figure 1 for the reference fine-scale solution.

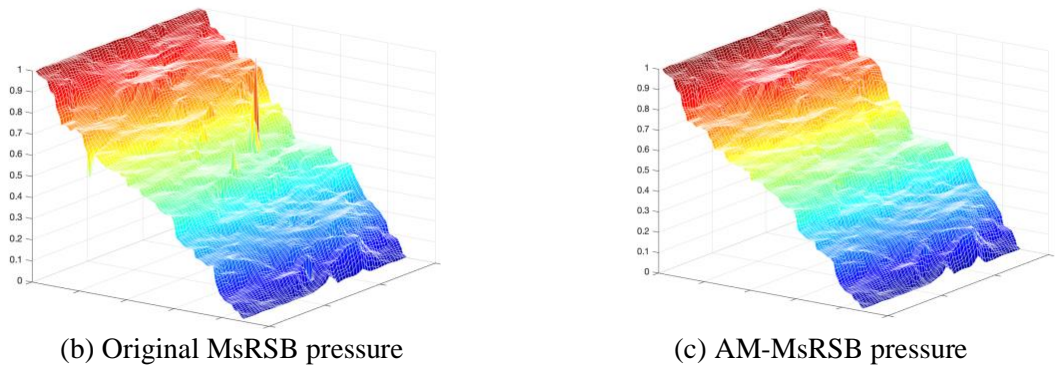

(b) Original MsRSB pressure        (c) AM-MsRSB pressure

**Figure 7.** Pressure surface plots for (a) the original MsRSB; and (b) the AM-MsRSB.

### 4.1.3. Test Case 3: 3D – full SPE10 problem

Next, we test the performance of AM-MsRSB on what has become a classical multiscale test case, the full 3D SPE 10 dataset. It consists of a total of around 1.1 million cells, specifically $(60 \times 220 \times 85)$ Cartesian grids. We use a coarsening ratio of $(5 \times 5 \times 5)$ and impose a unit pressure to drive flow from the left to the right boundary. **Figure 8** shows the reference fine-scale solution, MsRSB approximate solution, AM-MsRSB approximate solution, and the regions where MsRSB shows out-of-bound values; values above one are shown in burgundy and values below zero are shown in blue. With $\varepsilon = 0.0001$ and $w = 1$, AM-MsRSB operator produces strictly monotone and within bound solutions that are in very good agreement with the reference fine-scale solution ($L^2 = 0.0505$; $L^\infty = 0.269$). Meanwhile, the original MsRSB shows non-monotonic behavior as multiple regions show non-physical pressure values.

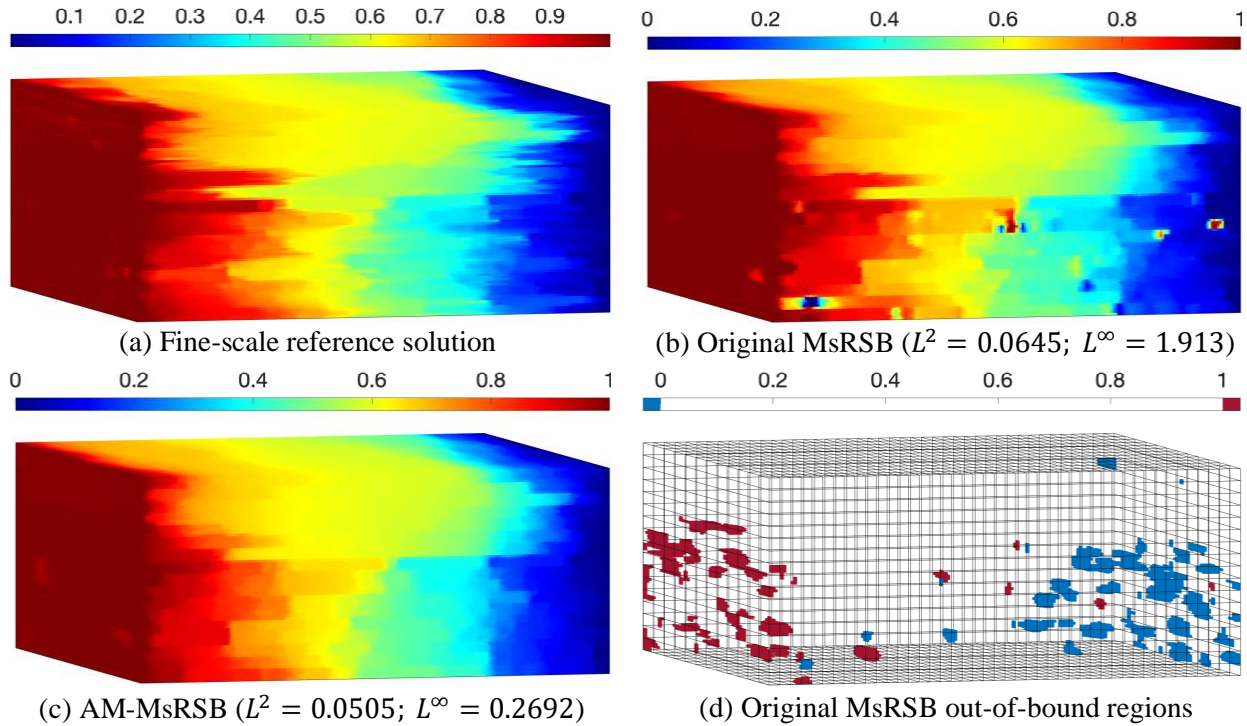

(a) Fine-scale reference solution  
(b) Original MsRSB ($L^2 = 0.0645$; $L^\infty = 1.913$)  
(c) AM-MsRSB ($L^2 = 0.0505$; $L^\infty = 0.2692$)  
(d) Original MsRSB out-of-bound regions

**Figure 8.** A comparison between (a) fine-scale reference solution, (b) MsRSB solution, (c) AM-MsRSB solution, and (d) regions where MsRSB produces out-of-bound values.

### 4.2. Multiscale methods as iterative solvers

The previous examples solely employed multiscale methods as approximate solvers. Nevertheless, as mentioned earlier, multiscale methods can also be used in an iterative setup to reduce the residual at the fine-scale. We investigate the performance of AM-MsRSB as an iterative solver on the 85[th] layer of the SPE10 problem, and we compare results with the original MsRSB. While it is always possible to achieve strict monotonicity with the AM-MsRSB operator, we note here that this may not always be the optimal choice within an iterative framework; Especially when the original MsRSB operator shows good monotonicity properties or when it barely violates the maximum principle. Therefore, although we will continue to refer to our algorithm as AM-MsRSB throughout this work, our selection for $\varepsilon$ and $w$ may not necessarily aim to achieve strictly monotone solutions.

It is often suggested to employ the multiscale operator as a two-stage preconditioner in a Krylov method like Generalized Minimal Residual Method (GMRES) to accelerate convergence. However, in this study, we will explore the performance of the AM-MsRSB as a stand-alone non-accelerated solver. As mentioned earlier, iterative multiscale formulations rely on an inexpensive relaxation method (smoother) for the updates. Herein, we will use incomplete LU-factorization with zero fill (ILU0). A multiscale cycle typically consists of one pass of the smoother and one multiscale coarse grid correction. Unless stated otherwise, we will be using one post-smoothing step.

### 4.2.1. Test case 4(a): SPE 10 bottom layer – iterative MsRSB

For this test case, we use the same problem setup of Test case 2 with the coarsening ratio of $(3 \times 5)$. Starting with an initial guess of $\boldsymbol{p}_f^0 = 0$, we aim to solve the system to a tight tolerance of $10^{-8}$. We compare the performance of MsRSB with AM-MsRSB for both restriction choices: the control volume $\boldsymbol{R}_{CV}$ and the Galerkin-type $\boldsymbol{R}_G$ operators.

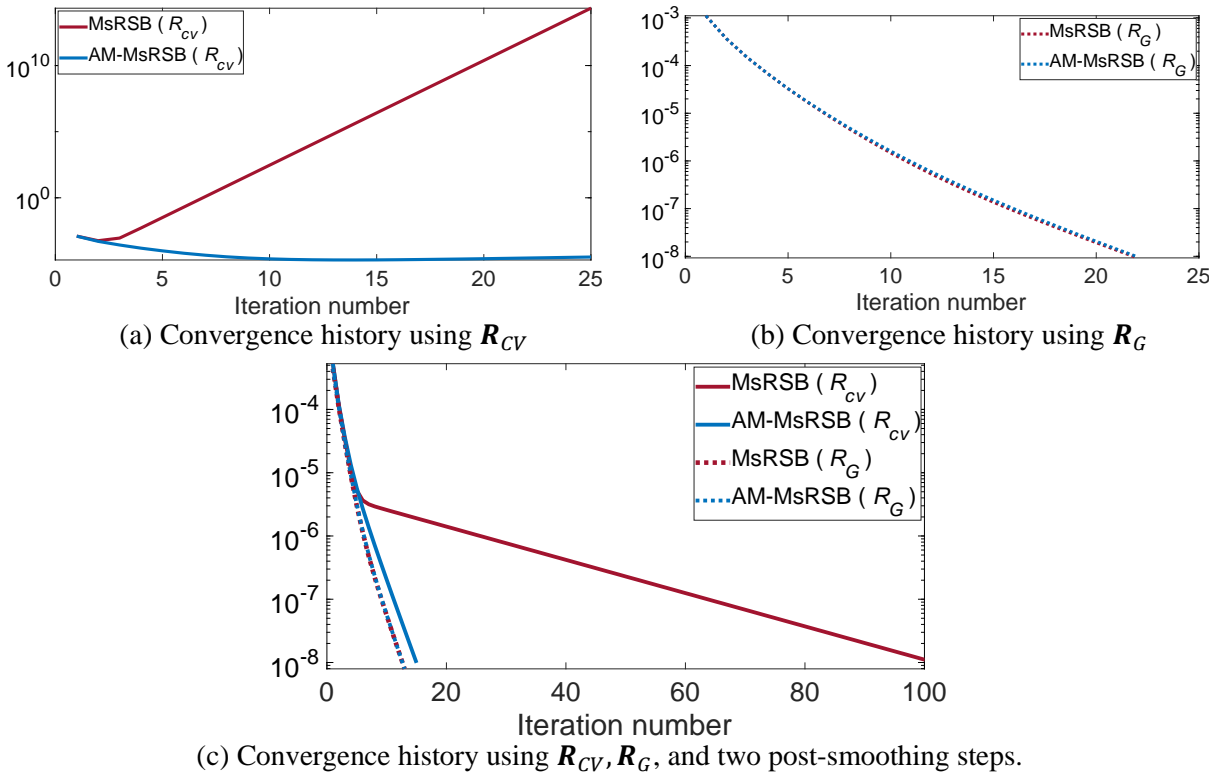

**Figure 9.** Convergence history for MsRSB and AM-MsRSB when used with (a) the mass-conservative control volume restriction operator $\boldsymbol{R}_{CV}$, (b) the Galerkin-type restriction operator $\boldsymbol{R}_G$ and (c) both restriction operators with two post-smoothing steps.

**Figure 9** reports the convergence performance of the two multiscale operators when used with $R_{CV}$, $R_G$, one post-smoothing step, and two post-smoothing steps. When $R_{CV}$ is employed, neither of the methods converge to a tolerance of $10^{-8}$. The original MsRSB operator diverges, while the AM-MsRSB stagnates around a residual of $10^{-5}$. **Figure 9b** illustrates the convergence of both operators when utilized with the Galerkin-type restriction operator. This choice aligns with literature recommendations, highlighting its advantage when pressure stability is a concern. When $R_G$ is used, the AM-MsRSB did not show an additional benefit compared to the original MsRSB and they both needed 22 iterations to converge. Next, we repeat the test with two post-smoothing steps, see **Figure 9c**. Results show that the AM-MsRSB reduced the fine-scale residual eight orders of magnitude within 15 iterations as opposed to 102 iterations using the original MsRSB with $R_{CV}$. When we replace $R_{CV}$ with $R_G$, both operators converge within 13 iterations.

### 4.2.2. Test case 4(b): SPE 10 bottom layer – iterative MsRSB

In this test, we investigate the effect of varying the coarsening ratio to $(7 \times 15)$, which leads to a reduced system of $(9 \times 15)$ coarse grids. A larger coarsening ratio implies that out-of-bound coarse-scale solutions, even for few coarse cells, can result in more violations when prolongated to the fine-scale system. As in the previous example, we start with an initial guess of $p_f^0 = 0$, and we target a tight tolerance of $10^{-8}$. A comparison of the convergence history is depicted in **Figure 10**.

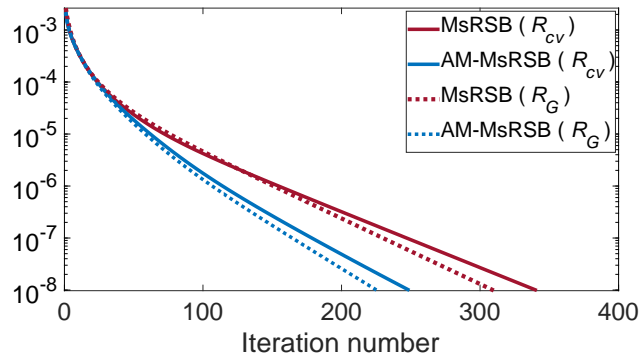

**Figure 10.** Comparison of the convergence history for test case 4(b) when using MsRSB and AM-MsRSB with (a) the mass-conservative control volume restriction operator $R_{CV}$ and (b) the Galerkin-type restriction operator $R_G$.

Like in the previous test case, 4(a), using the Galerkin-type restriction operator $R_G$ is advantageous in an iterative multiscale setup. The AM-MsRSB requires significantly less iterations to reduce the fine-scale residual to $10^{-8}$ when used with the control volume restriction operator. Furthermore, when used with the Galerkin-type restriction operator, AM-MsRSB converged within 226 iterations, compared to approximately 310 iterations for the original MsRSB. The high number of iterations is attributed to the large coarsening ratio which hinders the performance of both operators. However, the AM-MsRSB operator shows significant enhancements when used with large coarsening ratios. As mentioned earlier, convergence can be accelerated if the multiscale operator is used as a two-stage preconditioner in a Krylov method like GMRES.

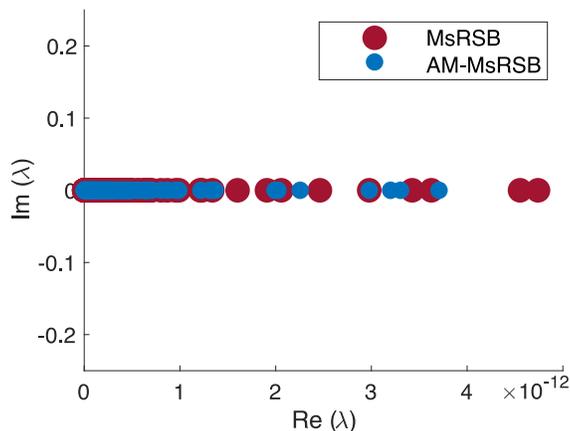

**Figure 11.** A zoomed-in screenshot of the lower end of the spectra of MsRSB and AM-MsRSB.

In **section 3**, we emphasized the significance of preserving the near null space of the original MsRSB operator, as it plays a crucial role in capturing the system's smooth modes. In **Figure 11,** we show a zoomed-in screenshot of the lower end of the spectra of MsRSB coarse-scale matrix, and AM-MsRSB coarse-scale matrix for this particular test case. It is evident that our algorithm effectively maintains the near null space, enabling the AM-MsRSB operator to accurately address low-frequency errors.

## 5. Multiscale methods for MPFA discretized fine-scale systems

The bulk of work on multiscale finite volume methods mainly targeted TPFA discretized fine-scale systems. This is because the TPFA scheme continues to be the default setting for many commercial and research simulators. Lately, however, there has been a growing interest in applying MsRSB to MPFA discretizations [45], [58], [59]. Based on the research conducted by Bosma et al., [58] it was found that the original MsRSB encounters problems when applied to non M-matrices, as it tends to generate divergent basis functions. To address this issue, the authors introduced an enhanced MsRSB approach, which involves a filtering strategy to enforce M-matrix properties on the fine-scale system prior to constructing the basis functions.

In this section, we modify our **Algorithm 1** slightly to address the issue of MsRSB with multipoint discretizations. As opposed to enforcing M-matrix properties at the fine-scale by dropping the off-diagonal positive entries, we propose redistributing the contribution of the positive entries in a similar manner to what we did earlier. So, we will introduce a perturbed fine-scale matrix $\widehat{A}_f = A_f + B$ with enforced M-matrix properties. Here, $B$ takes the exact same form as shown in equation **(12)** but with fine-scale indices. **Algorithm 2** summarizes our proposed methodology for enforcing M-matrix properties at the fine-scale.

---
**Algorithm 2** Obtaining improved basis functions for MPFA discretized systems
---
1: Locate all positive off-diagonal entries of $A_f$, $Pos\_entry$
2: **for** $i = 1: numel(Pos\_entry)$ **do**
3:     Construct the perturbation matrix $B$ as follows:
4:     $b_{ij} \leftarrow b_{ij} - a_{ij}^f$
5:     $b_{ii} \leftarrow b_{ii} + a_{ij}^f$
6:     $b_{ji} \leftarrow b_{ji} - a_{ij}^f$
7:     $b_{jj} \leftarrow b_{jj} + a_{ij}^f$
8: **end for**
9: $\widehat{A}_f = A_f + B$
10: Compute restricted-smoothed basis functions on $\widehat{A}_f$
---

The perturbed fine-scale matrix $\widehat{A}_f$ is only used for computing the basis functions. When constructing the reduced coarse-scale system, we use the intact fine-scale matrix, $A_f$. Using our algorithm to enforce M-matrix properties at the fine-scale generates partition-of-unity prolongation operator with in-bound values [0, 1], which coincides with the findings of [58]. Note that **Algorithm 2** is meant to provide accurate and convergent basis functions, but when applied alone, it does not guarantee strictly monotone solutions. When monotone solutions are desired, **Algorithm 1** should be invoked at the coarse level.

### 5.1. Test Case 5(a): MsRSB for MPFA discretized system

To assess our algorithm's performance, we examine a curvilinear grid by perturbing the interior nodes of a $100 \times 100$ structured grid, which discretizes a 20 by 150 meters domain. The coarsening ratio applied is $(5 \times 5)$, resulting in $(20 \times 20)$ coarse grids in the reduced system. The permeability tensor $K$ has diagonal values of 100 mD and off-diagonal values of 75 mD. The fine-scale system is discretized using MPFA-O. Next we compute the basis functions using the original MsRSB method for three different approaches: (a) when the fine-scale matrix remains unchanged; (b) using enhanced MsRSB i.e., when M-matrix properties are imposed on the fine-scale matrix using the filtering strategy presented in [58]; and (c) when we apply our algorithm to enforce M-matrix properties by redistributing the positive entries' contribution to the diagonal entries (AM-MSRSB). **Figure 12** shows the computed basis functions for an internal cell for the three approaches after 5 weighted Jacobi iterations. To further test the performance, we evaluate the convergence characteristics of the multiscale operators within an iterative formulation. **Figure 13** shows the convergence history for the approaches.

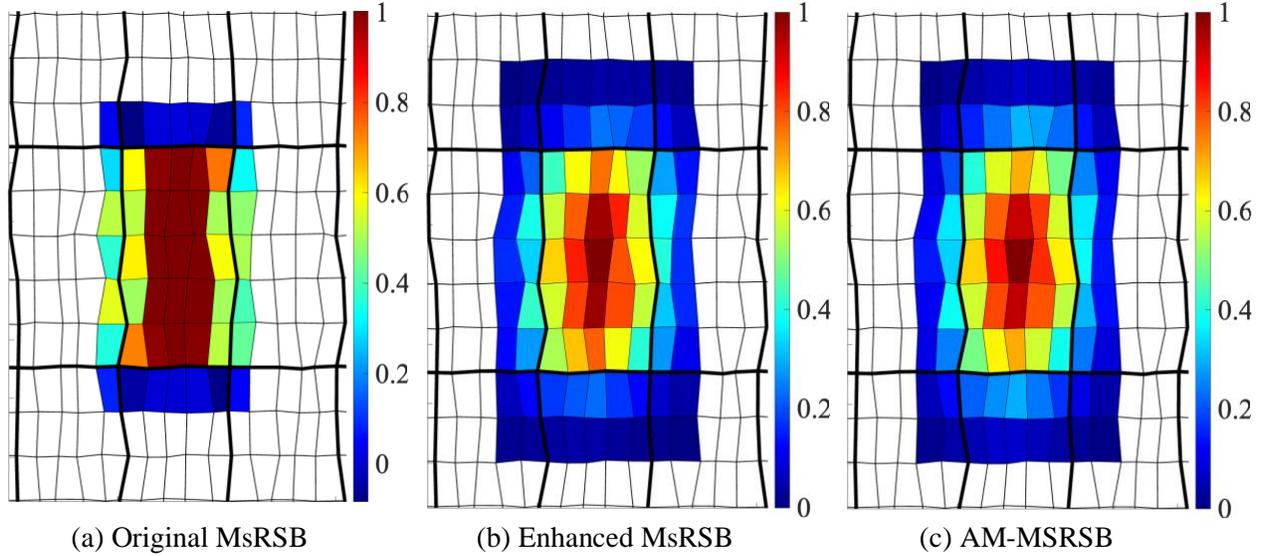

(a) Original MsRSB        (b) Enhanced MsRSB        (c) AM-MSRSB

**Figure 12.** Computed basis functions after 5 iterations for three different approaches: (a) the original MsRSB; (b) when filtering the fine-scale matrix; and (c) when applying our algorithm to enforce M-matrix properties on the fine-scale matrix.

Constructing basis functions on a fine-scale matrix with enforced M-matrix properties results in precise and bounded basis functions that significantly enhance the convergence speed of the multiscale solvers. In comparison to the filtering method, our approach not only matches but also slightly surpasses its performance.

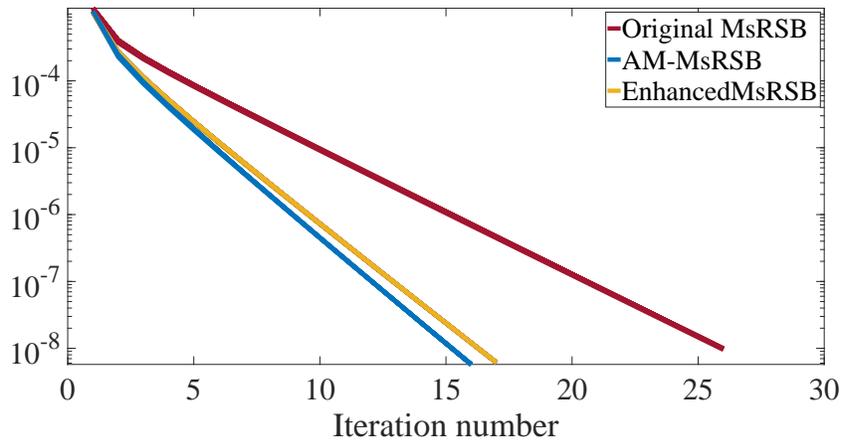

**Figure 13.** Comparison of the iterative performance for the three approaches.

## 5.2. Test Case 5(b): MsRSB for MPFA discretized system

In this test case, we introduce a rough curvilinear grid by perturbing all interior nodes of a $100 \times 100$ structured grid that discretizes a rectangular domain of 200 by 20 meters. We use a coarsening ratio of $(10 \times 10)$. **Figure 14** illustrates the curvilinear grid used for this case.

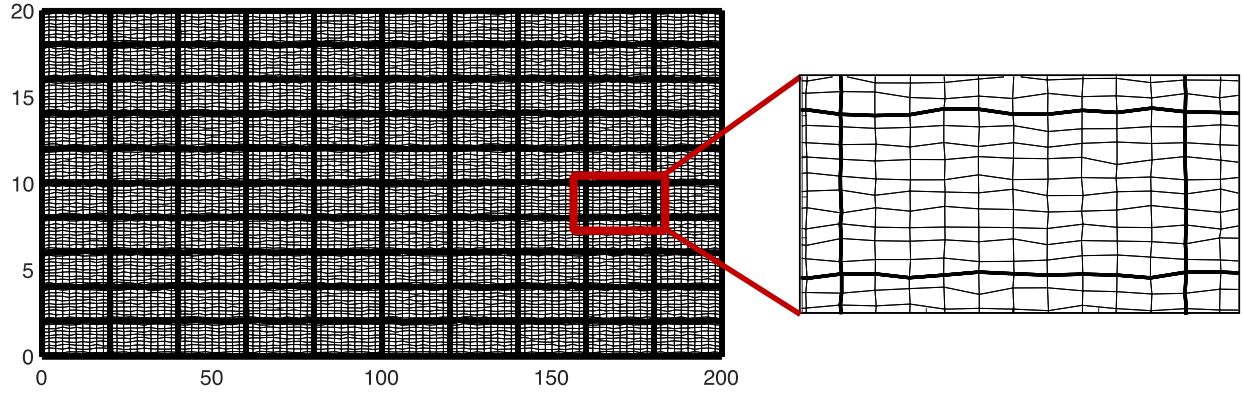

**Figure 14.** A curvilinear grid created by perturbing the interior nodes of regular $100 \times 100$ Cartesian grid.

The permeability is homogeneous but anisotropic and it is given by

$$K = \begin{bmatrix} cos\theta & -sin\theta \\ sin\theta & cos\theta \end{bmatrix} \begin{bmatrix} 1000 & 0 \\ 0 & 100 \end{bmatrix} \begin{bmatrix} cos\theta & sin\theta \\ -sin\theta & cos\theta \end{bmatrix},$$

where $\theta = 60°$. **Figure 15** shows the computed basis functions for an internal cell for the three approaches after 3, 5, and 10 weighted Jacobi iterations.

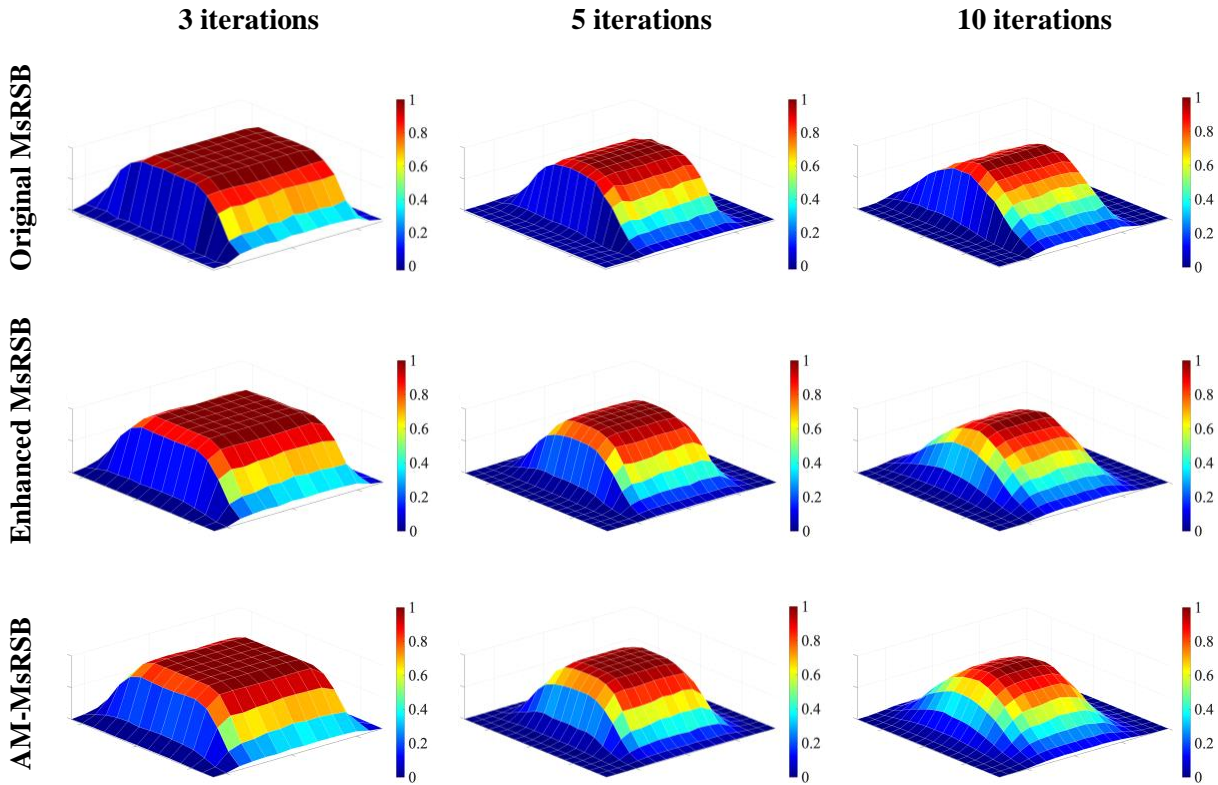

**Figure 15.** Computed basis functions after 3, 5, and 10 iterations for three different approaches: (a) the original MsRSB; (b) enhanced MsRSB; and (c) AM-MsRSB.

Although tough to eyeball, the original MsRSB produces basis functions that go out-of-bounds, whereas with the two other approaches all basis functions are within bounds. Moreover, treating the fine-scale matrix with our approach provides enhanced basis functions compared to the two other approaches. The basis functions obtained with our approach adjust quicker to the features of the underlying fine-scale flow equations.

To confirm this, we test the convergence properties of the multiscale operators in an iterative formulation. Starting with an initial guess of $\boldsymbol{p}_f^0 = 0$, we aim to reduce the fine-scale tolerance to $10^{-8}$. We use a single ILU0 post-smoothing step. **Figure 16** illustrates the convergence history for the three approaches. With our approach the residual is reduced by eight orders of magnitude within 54 iterations as opposed to 73 and 57 iterations for the original MsRSB and MsRSB with filtering, respectively.

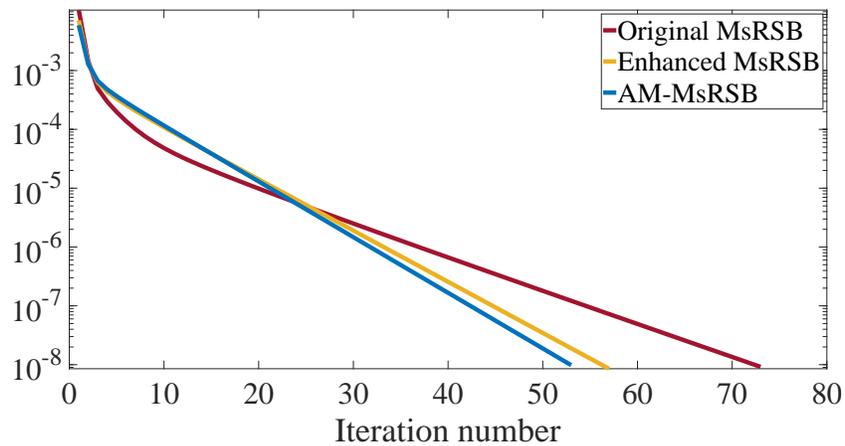

**Figure 16.** Convergence history for the three approaches.

### 5.3. Test case 5(c): MsRSB for MPFA discretized system

Next, we introduce a harsher problem with another $100 \times 100$ curvilinear grid discretizing a 500 by 200 meter domain. We choose a reduced system of 100 coarse cells, each containing 10 fine cells in each direction. A strong permeability anisotropy ratio of 100 is used ($\theta = 45°$):

$$\boldsymbol{K} = \begin{bmatrix} cos\theta & -sin\theta \\ sin\theta & cos\theta \end{bmatrix} \begin{bmatrix} 1000 & 0 \\ 0 & 10 \end{bmatrix} \begin{bmatrix} cos\theta & sin\theta \\ -sin\theta & cos\theta \end{bmatrix}.$$

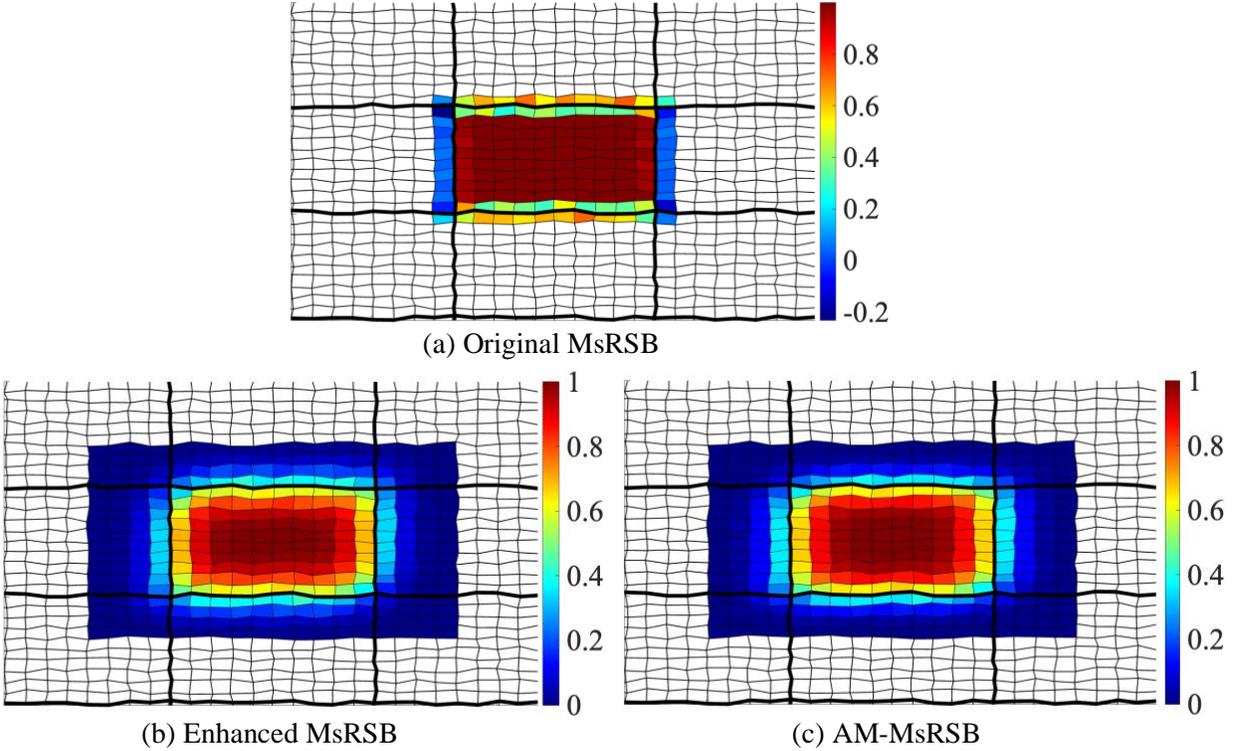

(a) Original MsRSB

(b) Enhanced MsRSB

(c) AM-MsRSB

**Figure 17.** Computed basis functions for an internal cell for the three different approaches.

In **Figure 17** we show the computed basis functions for each approach after 5 iterations and in **Figure 18** we report the convergence history of the multiscale operator for three different approaches of constructing the basis functions. The AM-MsRSB converges in 128 iterations as opposed to 143 iterations for the enhanced MsRSB operator. Both of which beat the original MsRSB which struggles to converge to $10^{-8}$.

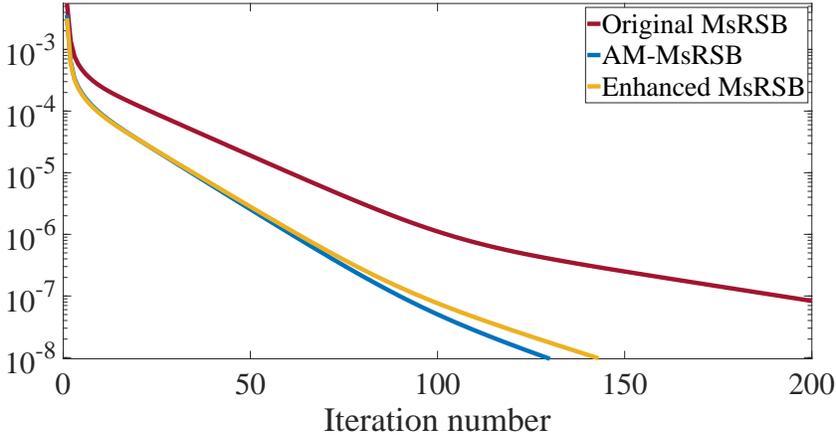

**Figure 18.** Convergence history for the multiscale operators for three different approaches of constructing the prolongation matrix.

## 5.4. Test Case 5(d): MsRSB for MPFA discretized system

In this final example, we test AM-MsRSB for MPFA discretized systems with the objective of obtaining monotone approximate solutions. As such, we first apply **Algorithm 2** to the fine-scale matrix to obtain accurate and convergent basis functions, then we invoke **Algorithm 1** to eliminate any non-physical oscillations that might arise from the multiscale operator. For the test case, the internal vertices of the SPE 10 bottom layer are perturbed. The $(60 \times 220)$ fine-scale grid is partitioned into $(12 \times 44)$ coarse cells. **Figure 19** shows the partitioned domain with permeability distribution in the background.

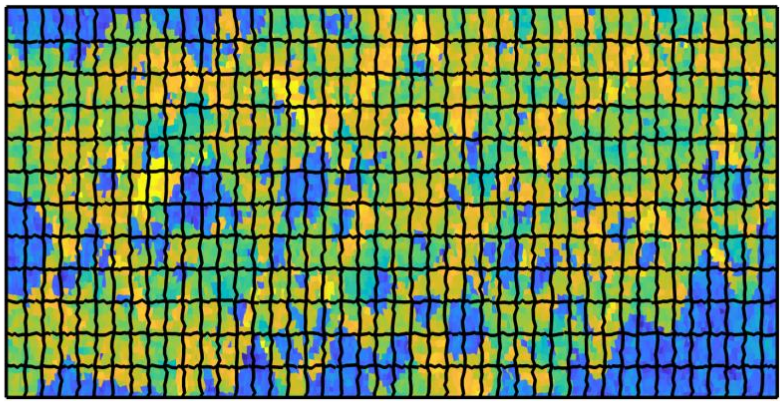

**Figure 19.** SPE 10 bottom layer with perturbed internal vertices.

The reduced system is then solved once at the coarse-scale and the solution is prolongated to the fine-scale. **Figure 20** shows the coarse-scale solutions obtained by the original MsRSB, the enhanced MsRSB i.e., with a filtered fine-scale system, and AM-MsRSB along with the prolongated pressure surfaces. The original MsRSB, overall, produced a smooth and an almost monotone solution except for the region near the outlet where it shows oscillatory behavior and goes out-of-bound. The enhanced MsRSB shows local non-physical oscillations and goes out-of-bound near the inlet boundary. In contrast, our AM-MsRSB provides a strictly monotone solution with very good accuracy at the fine-scale compared to the reference solution.

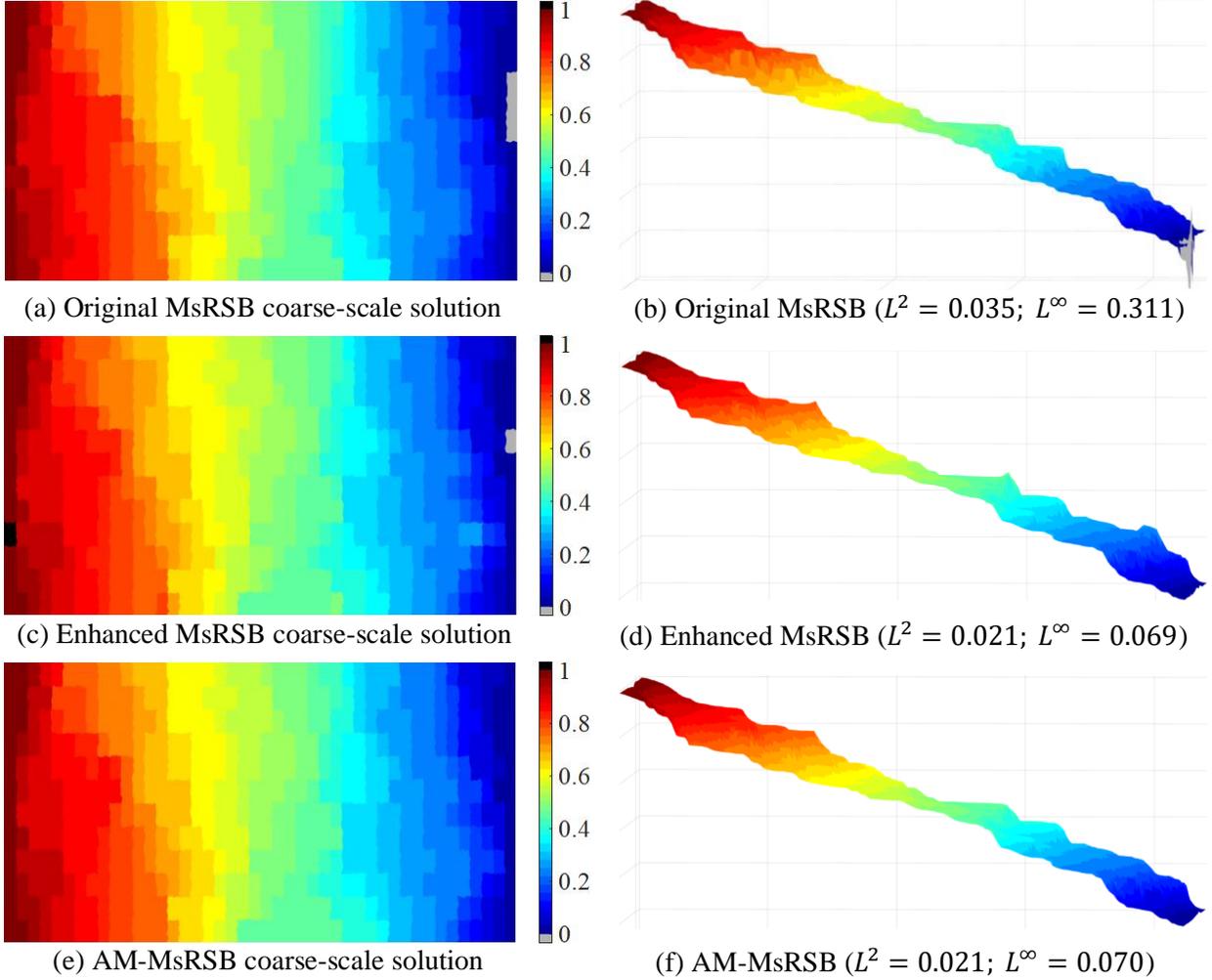

(a) Original MsRSB coarse-scale solution

(b) Original MsRSB ($L^2 = 0.035$; $L^\infty = 0.311$)

(c) Enhanced MsRSB coarse-scale solution

(d) Enhanced MsRSB ($L^2 = 0.021$; $L^\infty = 0.069$)

(e) AM-MsRSB coarse-scale solution

(f) AM-MsRSB ($L^2 = 0.021$; $L^\infty = 0.070$)

**Figure 20.** A comparison of the coarse-scale and fine-scale approximate solutions computed by the original MsRSB, the enhanced MsRSB, and AM-MsRSB.

## 6. Conclusions

In this work, a novel solution to the monotonicity issue of multiscale methods was proposed. Our algorithmic approach demonstrated significant improvements across all three applications of multiscale methods and for both TPFA and MPFA discretized systems.

The AM-MsFV/MsRSB guarantees monotonic and within bounds approximate solutions without compromising accuracy. Hence, our method is an ideal choice for multiscale methods as a one-step approximate solver. Furthermore, AM-MsFV/MsRSB yields monotone and precise solutions across different coarsening ratios, effectively addressing the challenge of multiscale methods' sensitivity to coarse grid partitioning choices. In iterative formulations where local conservation is desired, our monotone operator can be used as a concluding step with $\boldsymbol{R}_{CV}$ to provide strictly monotone fine-scale pressure solutions. By preserving the near null space of the original operator, the AM-MsRSB showed excellent performance when integrated in iterative formulations using both the control volume and the Galerkin-type restriction operators.

We further examined the AM-MsRSB performance on MPFA discretized systems. Our method proved to be quite potent and provides monotone and accurate approximate solutions when compared to the fine scale MPFA discretized solutions. Additionally, we explored the application of our approach to enhance the performance of the original MsRSB on MPFA discretized fine-scale systems, particularly in the construction of the prolongation operator. With AM-MsRSB, basis functions quickly adapt to the underlying fine-scale flow properties, which can result in computational cost savings during prolongation operator construction. In iterative setups, our approach reduced the residual to a very tight tolerance with fewer iterations compared to the original MsRSB and the enhanced MsRSB proposed in the existing literature. Finally, because our approach is algorithmic at the algebraic level, it is quite simple to implement in existing multiscale finite volume frameworks.